\documentclass[notitlepage,leqno,10pt]{article}

\usepackage[latin1]{inputenc}
\usepackage{amsfonts}
\usepackage{amsmath}
\usepackage{amssymb}
\usepackage{amsrefs}
\usepackage{color}
\definecolor{green2}{rgb}{0,0.6,0}

\newcommand{\R}{\mathbb{R}}
\newcommand{\N}{\mathbb{N}}
\newcommand{\proof}[1]{\par\smallskip\noindent{\bf Proof#1.}}
\newcommand{\qed}{\penalty 500\hfill$\square$\par\medskip}

\newcommand{\eps}{\varepsilon}

\def\XXint#1#2#3{{\setbox0=\hbox{$#1{#2#3}{\int}$}
     \vcenter{\hbox{$#2#3$}}\kern-.5\wd0}}

\def\sub{\underline}
\def\super{\overline}
\newtheorem{lem}         {Lemma}[section]
\newtheorem{conj}    [lem]{Conjecture}
\newtheorem{pro}    [lem]{Proposition}
\newtheorem{thm}    [lem]{Theorem}
\newtheorem{rem}    [lem]{Remark}
\newtheorem{cor}    [lem]{Corollary}

\title{Boundary blow-up solutions in the unit ball : asymptotics, uniqueness and symmetry}
\date{\today}
\author{O. Costin$^{1}$ and L. Dupaigne$^{2}$}
\begin{document}
\maketitle

{\begin{center}
$^{1}${\small
Department of Mathematics, The Ohio State University,\\
100 Math Tower, 231 West 18th Avenue, Columbus, OH 43210-1174, USA\\ {\tt costin@math.ohio-state.edu}\\
\smallskip
$^2$LAMFA, UMR CNRS 6140, Universit\'e Picardie Jules Verne \\33, rue St Leu, 80039 Amiens, France \\  \tt louis.dupaigne@math.cnrs.fr\\
}
\end{center}}
%%%%%%%%%%%%%%%%%%%%%%%%%%%%%%%%%%%%%%%%%%%%%%%%%%%%%%%%%%%%%%%%

\abstract{We calculate the full asymptotic expansion of boundary blow-up solutions (see equation \eqref{main radial} below), for any nonlinearity $f$. Our approach enables us to state sharp qualitative results regarding uniqueness and radial symmetry of solutions, as well as a characterization of nonlinearities for which the blow-up rate is universal. Lastly, we study in more detail the standard nonlinearities $f(u)=u^p$, $p>1$.}

\section{Introduction}
%Consider $\Omega$ a smoothly bounded domain of $\R^N$, $N\ge1$ and a 
Let $B$ denote the unit ball of  $\R^N$, $N\ge1$ and let $f\in C(\R)$. We study the equation
\begin{equation}\label{main radial} 
 \left\{
 \begin{aligned} 
\Delta u &=f(u) &\quad\text{in $B$,}\\
u &= +\infty &\quad\text{on $\partial B$,}
\end{aligned}
\right. 
 \end{equation} 
where the boundary condition is understood in the sense that
$$
\lim_{x\to x_{0}, x\in B}u(x) = +\infty\qquad\text{for all $x_{0}\in\partial B$}
$$
and where $f$ is assumed to be positive at infinity, in the sense that
\begin{equation} \label{positivity} 
\exists\; a\in\R\quad\text{s.t.}\quad f(a)>0\quad\text{and} \quad f(t)\ge0\quad\text{for $t>a$.}
\end{equation} 
A function $u$ satisfying \eqref{main radial} is called a boundary blow-up solution or simply a large solution. Existence of a solution of \eqref{main radial} is equivalent to the so-called Keller-Osserman condition : 
\begin{equation}\label{KO}
\int^{+\infty}\frac{dt}{\sqrt{F(t)}}<+\infty,\qquad\text{where $F(t)=\int_{a}^{t}f(s)\;ds$.}
\end{equation} 
For a proof of this fact, see the seminal works of J.B. Keller \cite{keller} and R. Osserman \cite{osserman} for the case of monotone $f$, as well as \cite{ddgr} for the general case. From here on, we always assume that \eqref{KO} holds. 

Our goal here is to study asymptotics, uniqueness and symmetry properties of solutions. 
%Furtherm	linearities $f$ according to the asymptotics they generate and provide relevant examples, as well as extensions of some of our results to smoothly bounded domains. 
Our approach improves known results in at least two directions : firstly, aside the necessary condition \eqref{KO}, we need not make \it any \rm additional assumption on $f$ to obtain the sharp asymptotics of solutions. Secondly, we obtain the complete asymptotic expansion of solutions, to \it all \rm orders.
Here is a summary of our findings.
\begin{thm} \label{theorem two} Let $f\in C(\R)$ and assume \eqref{positivity}, \eqref{KO} hold. Consider two solutions $u_{1}, u_{2}$  of \eqref{main radial}. Then,
$$
\lim_{x\to x_{0},x\in B}{u_{1}(x) - u_{2}(x)} = 0, \qquad\text{for all $x_{0}\in\partial B$.}
$$
More precisely, there exists a constant $C=C(u_{1},u_{2},N,F)>0$, such that for all $x\in B$,
\begin{equation} \label{first estimate} 
\left| u_{1}(x) - u_{2}(x) \right| \le C \int_{u_{2}(x)}^{+\infty}\frac{dt}{F(t)}\;dt.
\end{equation}   
In addition,
\begin{equation} \label{second estimate}
\left|  F(u_{1}) - F(u_{2}) \right| \in L^\infty(B).
\end{equation}   
\end{thm}
Estimates on the gradient of solutions can be obtained for a restricted class of nonlinearities, namely
\begin{thm}\label{theorem finer asymptotics} Let $f\in C(\R)$ and assume \eqref{positivity}, \eqref{KO} hold. Assume in addition that $f$ is increasing up to a linear perturbation i.e. there exists an increasing function $\tilde f$ and a constant $K$ such that
\begin{equation}\label{assumption gnn}
f(t) = \tilde f(t) - Kt, 
\end{equation}   
for all $t\in\R$. %in the range of $u$.
Consider two solutions $u_{1}, u_{2}$  of \eqref{main radial}. Then,
\begin{equation} \label{third estimate}
\left| \nabla (u_{1}-u_{2}) \right| \in L^\infty(B).
\end{equation} 
\end{thm}
%\begin{rem}
%Note that \eqref{KO2} is implied by \eqref{assumption gnn}.  
%\end{rem}
As a direct consequence of Theorem \ref{theorem two}, we obtain the following uniqueness result.
\begin{cor}\label{corollary one}Let $f\in C(\R)$ and assume \eqref{positivity}, \eqref{KO} hold. Assume in addition that $f$ is nondecreasing. Then, there exists a unique large solution of \eqref{main radial}.  
\end{cor}
%\begin{rem}
%Note that \eqref{KO2} automatically holds when $f$ is nondecreasing. 
%\end{rem}
\begin{rem}Many uniqueness theorems have been established in the literature (see e.g. the survey \cite{bandle}), and they hold for a  general class of bounded domains $\Omega$. However, in all of these results, additional assumptions on $f$ are needed, such as convexity.
\end{rem}
\proof{ of Corollary \ref{corollary one}  }
Let $u_{1},u_{2}$ denote two large solutions. It suffices to prove that $u_{1}\le u_{2}$. Assume this is not the case and let $\omega=\{x\in B\; : \; w(x)>0\}\neq\emptyset$, where $w=u_{1}-u_{2}$. Working if necessary on a connected component of $\omega$, we may always assume that $\omega$ is connected. Using Theorem \ref{theorem two}, we see that $w$ solves
the equation
 \begin{equation*}
 \left\{
 \begin{aligned} 
\Delta w &=f(u_{1})-f(u_{2})\ge 0 &\quad\text{in $\omega$,}\\
w &= 0 &\quad\text{on $\partial\omega$.}
\end{aligned}
\right. 
 \end{equation*}  
By the Maximum Principle, $w\le 0$ in $\omega$, a contradiction. 
\hfill\qed
When $f$ is not increasing, uniqueness fails in general. One may ask however whether all solutions of \eqref{main radial} are radial. H. Brezis made the following conjecture.
\begin{conj}[\cite{brezis}]\label{brezis} 
Let $f\in C^1(\R)$ denote a function such that \eqref{positivity}, \eqref{KO} hold. Then, every solution of \eqref{main radial} is radially symmetric.  
\end{conj}
A. Porretta and L. V\'eron proved that this is indeed the case under the additional assumption that $f$ is asymptotically convex (see \cite{pv}). We improve their result as follows.
\begin{cor}\label{corollary two} 
Let $f\in C^1(\R)$ and assume \eqref{positivity}, \eqref{KO} hold. Let $u$ denote a solution of \eqref{main radial}. Assume in addition that, up to a linear perturbation, $f$ is increasing  (i.e. \eqref{assumption gnn} holds for some nondecreasing function $\tilde f$ and some constant $K$).
Then, $u$ is radially symmetric. Furthermore, $\frac{\partial u}{\partial r}>0$ in $B\setminus\{0\}$.
\end{cor}
\begin{rem}
In the setting of the classical symmetry result of B. Gidas, W. M. Ni and L. Nirenberg (see \cite{gnn}),  \eqref{assumption gnn} is also assumed in order to prove symmetry. In the same article, the authors give an example of a nonlinearity $f$ failing \eqref{assumption gnn} for which there do exist nonradial solutions of the equation. In the context of large solutions, we do not have such a counter-example. In fact, we expect that none exists i.e. we believe that Conjecture \ref{brezis}  holds true. But at this stage, we do not even know whether radial symmetry remains true for simple nonlinearities such as $f(u)=u^2(1+\sin u)$.  
\end{rem}
Corollary \ref{corollary two}  is a direct consequence of the moving plane method and Theorem \ref{theorem finer asymptotics}  :
\proof{ of Corollary \ref{corollary two} } Let $U$ denote a radial solution of \eqref{main radial}. It follows from \eqref{positivity} that $U$ is a nondecreasing function of $r= \left|  x \right| $ for $r$ close to $1^-$ and $\frac{dU}{dr}(r)\to+\infty$ as $r\to1^-$. By \eqref{third estimate}, we conclude that any solution $u$ of \eqref{main radial} satisfies $\frac{\partial u}{\partial r}(x)\to+\infty$ as $x\to\partial B$, while  the tangential part of the gradient of $u$ remains bounded. We then apply Theorem 2.1 in \cite{pv}.
\hfill\qed
In addition to the \it relative \rm asymptotic information given by \eqref{first estimate}, \eqref{second estimate}  and \eqref{third estimate}, the \it exact \rm asymptotic expansion of a solution can be calculated to all orders, as follows.
{
\begin{thm} \label{theorem asymptotics} 
Let $f\in C(\R)$ and assume \eqref{positivity}, \eqref{KO} hold. Let $U_{0}\in\R$, $I=[U_{0},+\infty)$ and let $v_{0}$ be the function defined for $u\in I$ by 
\begin{equation} \label{def v nought}
v_{0}(u) = \sqrt{2F(u)}.
\end{equation}  
Consider the Banach space 
$$
\mathcal X = \{ v\in C(I;\R)\; : \; \exists M>0 \text{ such that } \left| v \right| \le M v_{0} \},
$$ 
endowed with the norm $\| v \| = \sup_{I} \left| v/{v_{0}}\right| $. If  the constant $U_{0}$ is chosen sufficiently large, there exists $\rho\in(0,1)$ such that there exists a unique solution $v\in\mathcal B(v_{0},\rho)\subset \mathcal X$ of the integral equation
\begin{equation} \label{equation v} 
v(u) = \sqrt{2\left(F(u)-(N-1)\int_{U_{0}}^{u} \frac{v}{r}\;dt\right)},\qquad u\in I,
\end{equation} 
where $r=r(u,v)$ is given for $u\in I$, $v\in \mathcal B(v_{0},\rho)$ by
\begin{equation}\label{r of v}  
r(u,v) = 1 - \int_{u}^{+\infty}\frac1v \;dt.
\end{equation} 
In addition, $v$ is the limit in $X$ of $\left(v_{k}\right)$ defined for $k=0$ by \eqref{def v nought} and for $k\ge1$, by
\begin{equation} \label{expression vk}
v_{k}(u) = \sqrt{2\left(F(u)-(N-1)\int_{U_{0}}^{u} \frac{v_{k-1}}{1 - \int_{t}^{+\infty}\frac1{v_{k-1}} \;ds}\;dt\right)}
\end{equation}   
and the sequence $v_{k}$ is asymptotic to $v$ i.e. as $u\to+\infty$,
$$
v_{k+1}(u)= v_{k}(u) + o(v_{k}(u))
$$
and given any $k\in\N$,
$$
v(u)= v_k(u) + O(v_{k+1}(u)-v_{k}(u)).
$$
Let now $u$ denote any solution of \eqref{main radial} and fix $r_{0}\in (0,1)$ such that $u(x)\ge U_{0}$ for $\left| x \right| \ge r_{0}$. For $k\ge 0$, define $u_{k}$ for $r\ge r_{0}$ as the unique solution of
\begin{equation}\label{expression uk}
\left\{
\begin{aligned}  
\frac{du_{k}}{dr} &= v_{k}(u_{k})\\
\lim_{r\to1^-}u_{k}(r) &= +\infty,
\end{aligned}
\right.
\end{equation}  
where $v_{k}$ is given by \eqref{expression vk}. 
Then, as $r\to1^-,$
%$$u_{k+1}(r)= u_{k}(r+o(1-r))$$
$$
\int_{u_{k}(r)}^{u_{k+1}(r)}\frac{du}{v_{0}}= o\left(\int_{u_{k}(r)}^{+\infty}\frac{du}{{v_{0}}}\right)
$$
and given any $k\in\N$, we have as $x\to\partial B$, 
\begin{equation} \label{asymptotics of u} 
\begin{aligned}
\int_{u_{k}(\left| x \right| )}^{u(x)}\frac{du}{v_{0}}= o\left(\int_{u_{k}(\left| x \right| )}^{+\infty}\frac{du}{{v_{0}}}\right)
%u(x)=u_{k}(\left| x \right| + o(1-\left| x \right|)).% + O(u_{k+1}(\left| x \right| )-u_{k}(\left| x \right| )).
\end{aligned}
\end{equation} 
\end{thm}
Theorem \ref{theorem asymptotics} enables one to calculate (implicitely) the asymptotic expansion of a solution term by term. But how many terms in this expansion are singular? This is what we discuss in our last set of results.

We begin with the simplest class of nonlinearities $f$, those for which only one term in the expansion is singular, namely the function $u_{0}$ defined by \eqref{def v nought} and \eqref{expression uk}. It turns out, as A.C. Lazer and P.J. Mc Kenna first demonstrated (see \cite{lm}), that in this case $u_{0}(1-d(x))$ is the only singular term in the asymptotics of any blow-up solution on any smoothly bounded domain $\Omega\subset\R^N$ and for any dimension $N\ge1$, where $d(x)$ denotes the distance of a point $x\in\Omega$ to the boundary of $\Omega$. In other words, the blow-up rate is universal. The question is now to determine for which nonlinearities $f$, this universal blow-up occurs. We characterize these nonlinearities as follows:
\begin{thm}\label{theorem universal} Let $\Omega\subset\R^N$ denote a bounded domain satisfying an inner and an outer sphere condition at each point of its boundary. Let $f\in C(\R)$, assume \eqref{positivity}, \eqref{KO} hold and consider the equation
 \begin{equation}\label{main} 
 \left\{
 \begin{aligned} 
\Delta u &=f(u) &\quad\text{in $\Omega$,}\\
u &= +\infty &\quad\text{on $\partial\Omega$.}
\end{aligned}
\right. 
 \end{equation} 
Assume
\begin{equation}\label{assumption one term}
\lim_{u\to+\infty} \sqrt {2F(u)}\int_{u}^{+\infty}\frac{\int_{0}^{t} \sqrt{2F}\;ds}{(2F)^{3/2}}\;dt =0.
\end{equation}   
Then, any solution of \eqref{main} satisfies
 \begin{equation}\label{one term only}  
 \lim_{x\to\partial\Omega} u(x) - u_{0}(1-d(x)) = 0,
 \end{equation} 
 where $d(x)={\rm dist} (x,\partial\Omega)$ and $u_{0}$ is defined by \eqref{def v nought}, \eqref{expression uk}. 
 
We also have the following partial converse statement : if 
\begin{equation}\label{assumption one term fails}
\liminf_{u\to+\infty} \sqrt {2F(u)}\int_{u}^{+\infty}\frac{\int_{0}^{t} \sqrt{2F}\;ds}{(2F)^{3/2}}\;dt >0,
\end{equation}  
then \eqref{one term only} always fails. 
\end{thm}
 \begin{rem}
To our knowledge, \eqref{assumption one term} improves upon all known conditions for \eqref{one term only} to hold (see in particular \cite{lm} and \cite{bandle-marcus}) . Despite its unappealing technical looks, \eqref{assumption one term}  uses information on the asymptotics of $F$ only (in particular, no direct information on $f$ is required). Nonlinearities such that $F(u)\sim e^u$ or $F(u)\sim u^p$, $p>4$ as $u\to+\infty$ qualify. For $F(u)\sim u^4$, \eqref{assumption one term fails} holds and so the conclusion \eqref{one term only} fails.      
\end{rem}

\begin{rem}
Condition \eqref{assumption one term} can be weakened to : 
$$
\lim_{r\to1^-} \sqrt {2F(u_{0})}\int_{u_{0}}^{+\infty}\frac{\int_{0}^{t} \sqrt{2F}\;ds}{(2F)^{3/2}}\;dt =0,
$$
where $u_{0}=u_{0}(r)$ is defined by \eqref{expression uk}. Similarly, \eqref{assumption one term fails} can be weakened to  
$$
\liminf_{r\to1^-} \sqrt {2F(u_{1})}\int_{u_{1}}^{+\infty}\frac{\int_{0}^{t} \sqrt{2F}\;ds}{(2F)^{3/2}}\;dt >0.
$$
\end{rem}
As an immediate corollary, we obtain uniqueness on general domains, whenever only one singular term appears:
\begin{cor} Assume \eqref{assumption one term}. If in addition, $f$ is nondecreasing, then the solution of \eqref{main} is unique. 
\end{cor}
\proof{} Simply repeat the proof of Corollary \ref{corollary one}.\hfill\qed  
More than one term can be present in the asymptotic expansion of $u$. Finding all the (singular) terms in this expansion is of staggering algebraic complexity. To illustrate this, we provide the first three terms (in implicit form).
\begin{pro}\label{first three terms} Let $u_{2}$ be defined by \eqref{expression uk} for $k=2$. Let also $R_{1}$, $R_{2}$, $R_{3}$ denote three real-valued functions defined for $U\in\R$ sufficiently large by
\begin{multline*}
R_{0}(U) = \int^{+\infty}_{U}\frac{du}{\sqrt{2F}},\qquad
R_{1}(U) =(N-1)\int^{+\infty}_{U}\frac{\int^{u}\sqrt{2F}\;dt}{(2F)^{3/2}}\;du,\\
R_2(U)  = (N-1)\times\\
\times\int^{+\infty}_{U}
\left(\begin{aligned}
-\int^{u}\left((N-1)\frac{\int^{t}\sqrt{2F}\;ds}{\sqrt{2F}}+\sqrt{2F}\int^{+\infty}_{u}\frac{ds}{\sqrt{2F}}\right)\;dt+\\
+\frac{5(N-1)}4\frac{\left(\int^{u}\sqrt{2F}\;dt\right)^2}{2F}
\end{aligned}
\right)
\;\frac{du}{(2F)^{3/2}}. 
\end{multline*}
Then, for all $r\in(0,1)$, $r$ close to $1$, there holds
$$
1-r = R_{0}(u_{2}(r))+R_{1}(u_{2}(r))+R_{2}(u_{2}(r))(1+o(1)).
$$
\end{pro}
%\begin{rem} Both the second and the third order terms
%depend on the curvature of the ball, as the factor $(N-1)$ indicates. 
%\end{rem}

For specific nonlinearities, it is possible to invert the above identity. This is what we do for $f(u)=u^p$, $p>1$:
\begin{pro}\label{utothep} 
Let $p>1$ and $f(u)=u^p$, for $u>0$. Then, the unique positive solution of \eqref{main radial} satisfies as $r\to1^-$,
$$
u = d^{-\frac2{p-1}}\sum_{k=0}^{[2/(p-1)]}a_{k}d^k + o(1),
$$ 
where $d(r)=1-r$ for $r\in(0,1)$, and where each $a_{k}\in\R$ depends on $N$ and $p$ only.
\end{pro}

\begin{rem}Note that the above result remains true for any nonlinearity $f$ such that, for some positive constant $c$, $F(u)\sim c u^{p+1}$ as $u\to+\infty$ (and any solution of the equation).
\end{rem}

\

\noindent{\bf Outline of the paper}

\noindent \begin{enumerate}
\item In the next section, we show that any solution $u$ of \eqref{main radial} can be squeezed between two radial solutions $U$ and $V$ i.e. the inequality $U\le u\le V$ holds throughout $B$. 

\item Thanks to this result, we need only find the asymptotics of \it radial \rm solutions to prove Theorem \ref{theorem two}.  This is what we do in Section \ref{section 3}. 

\item To obtain gradient estimates, the squeezing technique is insufficient and more work is needed. In Section \ref{section 4}, we estimate tangential derivatives via a standard comparison argument, while we gain control over the radial component through a more delicate potential theoretic argument.   

\item Section \ref{section 5} is dedicated to the proof of Theorem \ref{theorem asymptotics}, that is we establish an algorithm for computing the asymptotics of  solutions to all orders.

\item In Section \ref{section 6}, we characterize nonlinearities for which the blow-up rate is universal.

\item At last, Sections \ref{section 7} and \ref{section 8} contain the tedious calculations of the first three terms of the asymptotic expansion of $u$ in implicit form for general $f$, and of all terms explicitely for $f(u)=u^p$.         
  \end{enumerate}

\noindent{\bf Notation}

\noindent Throughout this paper, the letter $C$ denotes a generic constant, the value of which is immaterial. In the last section of the paper, we use the symbol $c_{k}$ to denote a quantity indexed by an integer $k$, thought of being ``constant for fixed $k$'', the value of which is again immaterial. 

}

\section{Ordering solutions}
In this section, we prove that any solution of the equation is bounded above and below by  radial blow-up solutions. To do so, we impose the following additional condition: $g(t):=f(-t)$ satisfies \eqref{positivity} and 
\begin{equation}\label{KO2}
\int^{+\infty}\frac{dt}{\sqrt{G(t)}}=+\infty, \qquad\text{where $G'(t)=f(-t)$.}
\end{equation} 
\begin{rem}\label{remark KO2} 
Note that \eqref{KO2}  is not restrictive. Indeed, if $u$ denotes a solution of \eqref{main radial} and $m=\min_{B}u$, then $u$ also solves \eqref{main radial} with nonlinearity $\tilde f$ defined for $u\in\R$ by
\begin{equation*}
\tilde f(u)=\left\{
\begin{aligned}
f(m)+(m-u)&\qquad\text{if $u<m$}\\
f(u)&\qquad\text{if $u\ge m$.}
\end{aligned}
\right.
\end{equation*}
Then, $\tilde f$ clearly satisfies \eqref{KO2}.    
\end{rem}
We now proceed through a series of three lemmas.
\begin{lem}\label{lemma KO2}  Assume \eqref{KO2} holds. For $M\in\R$ sufficiently large, there exists a radial function $\sub v\in C^2(B)\cap C(\super B)$ satisfying
$$
\Delta \sub v \ge f(\sub v)\qquad\text{in $B$}
$$
and such that 
$$\sub v\le -M\qquad\text{in $B$.}$$
\end{lem}
 \proof{} Let $g(t)=f(-t)$ for $t\in\R$ and let $a>M$ a parameter to be fixed later on. Since $g$ satisfies \eqref{positivity}, we may always assume that $g(t)\ge0$ for $t\ge M$. Let now $w$ denote a solution of 
\begin{equation}\label{diff eq} 
 \left\{
 \begin{aligned} 
-w'' &=g(w) \\
w(0) &= a \\
w'(0) &=0
\end{aligned}
\right. 
 \end{equation}   
 
{\noindent \bf Claim. } There exists $a>M$ sufficently large, such that $w(1)\ge M$.

Note that $w$ is nonincreasing in the set $\{t\;:\; w(t)\ge M\}$. We distinguish two cases. 

{\noindent \bf Case 1.} $w>M$. 

In this case, $w$ is defined on all of $\R^+$. In particular, $w(1)> M$, as desired.

{\noindent \bf Case 2.} There exists $R>0$ such that $w(R)=M$. 

In this case, since $w$ is nonincreasing in $(0,R)$, we just need to prove that $R\ge1$. To do so, multiply \eqref{diff eq} by $-w'$ and integrate between $0$ and $r\in(0,R)$:
$$
-w' = \sqrt{2(G(a)-G(w))},
$$ 
where $G$ is an antiderivative of $g$. Integrate again between $0$ and $R$:
$$
\int_{M}^{a} \frac{dt}{\sqrt{2(G(a)-G(t))}}=\int_{0}^R \frac{-w'}{\sqrt{2(G(a)-G(w))}}\;dr = R.
$$
Now, 
\begin{multline*} 
R=\int_{M}^{a} \frac{dt}{\sqrt{2(G(a)-G(t))}}\ge \int_{M}^{G^{-1}(G(a)/2)} \frac{dt}{\sqrt{2(G(a)-G(t))}}\ge\\ \ge\int_{M}^{G^{-1}(G(a)/2)} \frac{dt}{\sqrt{2G(t)}}.
\end{multline*}
By \eqref{KO2}, we deduce that $R\ge1$ for $a$ sufficently large. We have just proved that $\left. w\right|_{(0,1)} \ge w(1)\ge M$ and the claim is proved.
 
It follows that the function $\sub v$ defined for $x\in B$ by $\sub v(x)=-w( \left| x \right| )$, is the desired subsolution.
\hfill\qed 
 
\begin{lem}\label{lemma KO} Let $f\in C(\R)$ and assume \eqref{positivity}, \eqref{KO} hold. Assume $\sub v~\in~%C^2(B)\cap
C(\super B)$ satisfies
$$
\Delta \sub v \ge f(\sub v)\qquad\text{in $B$.}
$$
Then, there exists a radial large solution $V$ of \eqref{main radial} such that $V\ge \sub v$.  
\end{lem}
\proof{}
Let $\super v:= N$. Then, $\sub v$ and $\super v$ are respectively a sub and supersolution of 
\begin{equation}\label{truncated problem}  
 \left\{
 \begin{aligned} 
\Delta v &=f(v) &\quad\text{in $B$,}\\
v &= N &\quad\text{on $\partial B$,}
\end{aligned}
\right. 
 \end{equation}     
 provided $N$ is chosen so large that $N> \| \sub v \|_{L^\infty(B)}$ and $f(N)\ge0$.
Futhermore, $\sub v< \super v$ in $B$ for such values of $N$. By the method of sub and supersolutions (see e.g. Proposition 2.1 in \cite{ddgr}), there exists a minimal solution $V_{N}$ of \eqref{truncated problem} such that $N\ge V_{N} \ge\sub v$.
 Note that $V_{N}$ is radial, as follows from the classical symmetry result of Gidas, Ni and Nirenberg (see \cite{gnn}) when $f$ is locally Lipschitz, or simply from the fact that $V_{N}$ is minimal, hence radial, for $f$ merely continuous. Also, since $V_{N}$ is minimal, we have  that the sequence $(V_{N})$ is nondecreasing with respect to $N$ (apply e.g. the Minimality Principle, Corollary 2.2, in \cite{ddgr}). 
 
 It turns out that the sequence $(V_{N})$ is uniformly bounded on compact sets of $B$. Indeed, fix $R_{1}<1$. There exists a solution $\tilde U$ blowing up on the boundary of the ball of radius $1$ and satisfying $\tilde U\ge \sub v$ in $B_{R_{1}}$, see Remark 2.9 in \cite{ddgr}.  
By minimality, $\sub v \le V_{N}\le \tilde U$ in $B_{R_{1}}$, whence $(V_{N})$ is uniformly bounded on $B_{R_{2}}$, for any given $R_{2}<R_{1}$. 

We have just proved that each $V_{N}$ is radial and that the sequence $(V_{N})$ is nondecreasing and bounded on compact subsets of $B$. By standard elliptic regularity, it follows that $(V_{N})$ converges to a radial solution $V$ of \eqref{main radial}, such that $V\ge\sub v$ in $B$. 
\hfill\qed
\begin{lem}\label{lemma three}  Assume \eqref{KO} and \eqref{KO2} hold. 
Let %$u\in C^2(B)$, 
$u$ denote a solution of \eqref{main radial}. 
Then, there exist two radial functions $U,V$ solving \eqref{main radial} such that 
$$
U\le u\le V \qquad\text{in $B$.}
$$ 
\end{lem}
\proof{} 
Let $-M$ denote the minimum value of $u$ and let $\sub v$ denote the subsolution given by Lemma \ref{lemma KO2}. In particular, $\sub v\le u$. By Lemma \ref{lemma KO}, there exists a solution $U\ge\sub v$ of \eqref{main radial} and we may asssume that $U$ is the minimal solution relative to $\sub v$ i.e.  given any other solution $\tilde u\ge\sub v$ of \eqref{main radial}, $U\le \tilde u$. In particular, $U \le u$. It remains to construct a radial solution $V$ of \eqref{main radial} such that $u\le V$. To do so, we fix $R<1$. By Lemma \ref{lemma KO}, letting $\sub v=\left.u \right|_{B_{R}}$, there exists a radial solution $v=V_{R}$ of 
\begin{equation}\label{problem br}  
 \left\{
 \begin{aligned} 
\Delta v &=f(v) &\quad\text{in $B_{R}$,}\\
v &= +\infty &\quad\text{on $\partial B_{R}$,}
\end{aligned}
\right. 
 \end{equation}  
such that $V_{R}\ge u$ in $B_{R}$. Since $V_{R}$ is constructed as the monotone limit of minimal solutions $V_{N}$ (see the proof of the previous lemma), one can easily check that the mapping $R\mapsto V_{R}$ is nonincreasing (hence automatically bounded on compact sets of $B$). Hence, as $R\to1$, $V_{R}$ converges to a solution $V$ of \eqref{main radial}, which is radial and  satisfies $V\ge u$ in $B$, as desired.
\hfill\qed
\section{Asymptotics of radial solutions}\label{section 3} 
Our next result establishes that the asymptotic expansion of a radial blow-up solution is unique. More precisely, consider the one-dimensional problem
\begin{equation} \label{1D} 
\frac{d^2 \phi}{dr^2} = f(\phi), \qquad r<1, \qquad \phi(r)\to+\infty\;\text{ as }\; r\to 1^-.
\end{equation} 
All solutions are given implicitely by
$$
\int_{\phi}^{+\infty}\frac{ds}{\sqrt {2F(s)}} = 1-r, \quad\text{ where }\quad F'=f. 
$$
We recall the following fact, first observed by C. Bandle and M. Marcus in \cite{bandle-marcus} : 
\begin{rem}Let $\phi$ and $\phi_{c}$ denote two solutions of \eqref{1D} corresponding to the antiderivatives $F$ and $F+c$, respectively. Then $\phi(r)-\phi_{c}(r)\to0$ as $r\to1^-$.
\end{rem}
We improve this result in the following way.
\begin{thm}\label{theorem one} 
Let $N\ge1$ and let $u_1, u_{2}$ denote two strictly increasing functions solving 
\begin{equation} \label{asymptotic 1} 
\left\{
\begin{aligned}
\frac{d^2 u}{dr^2}  + \frac{N-1}{r}\frac{du}{dr} &= f(u),\qquad\text{$r<1$,}\\
\lim_{r\to1^-} u(r) = +\infty.
\end{aligned}
\right.
\end{equation} 
Then, 
\begin{equation*}
\left| u_{1}(r) - u_{2}(r) \right| \le C \int_{u_{2}(r)}^{+\infty}\frac{dt}{F(t)}\;dt.
\end{equation*}   
In addition, the quantity
$\left|  F(u_{1}) - F(u_{2}) \right|$ is bounded.
\end{thm} 
\begin{rem}
Clearly, Theorem \ref{theorem two} follows as a direct consequence of Remark \ref{remark KO2},  Lemma \ref{lemma three} and  Theorem  \ref{theorem one}.  
\end{rem}

\proof{} We want to think of the second term on the left-hand side of equation \eqref{asymptotic 1} as a lower order perturbation as $r\to1$.  So, we integrate \eqref{asymptotic 1} in the same way we would solve \eqref{1D}, namely we let $v=du/dr$ and multiply the equation by $v$. We get
$$
\frac{d}{dr}\left(\frac{v^2}2 \right) + \frac{N-1}{r}v^2 = \frac{d}{dr}\left(F(u)\right).
$$
We define the resulting error term  by
\begin{equation} \label{definition g}
g : = -\frac{v^2}{2} + F(u),
\end{equation}   
which, seen as function of $r$, satisfies the differential equation
\begin{equation} \label{equation g in r}
\frac{dg}{dr} = \frac{N-1}{r}v^2.
\end{equation}  
Since $u$ is a strictly increasing function, the change of independent variable $u=u(r)$ is valid. Thinking of $g$ as a function of the variable $u$, we have $\frac{dg}{du}=\frac{dg}{dr}\frac{dr}{du}=\frac1{v}\frac{dg}{dr}$ and so
\begin{equation} \label{equation g in u}
\frac{dg}{du} =  \frac{N-1}{r}v.
\end{equation}  
Since \eqref{asymptotic 1} holds for $r$ close to $1$, the above equation holds for $u$ in a neighborhood of $+\infty$. Solving \eqref{definition g} for $v$, we finally obtain
\begin{equation}\label{equation g in u full}
\left\{
\begin{aligned}
\frac{dg}{du} &=  \frac{N-1}{r}v=\frac{N-1}{r}\sqrt{2(F(u)-g)},\\
\frac{dr}{du}&=1/v=\left( 2(F(u)-g)\right)^{-1/2}.
\end{aligned}
\right.
\end{equation}    
We start by calculating the leading asymptotic behaviour of $g$ at $+\infty$ :
\begin{lem} \label{lemma one}
$$
\lim_{u\to+\infty} \frac{g(u)}{F(u)} = 0.
$$  
In addition,
\begin{equation} \label{g is equivalent to G}  
\lim_{u\to+\infty}\frac{g(u)}{(N-1)G(u)}=1,\quad\text{where $G$ is any antiderivative of $\sqrt{2F}$.}
\end{equation} 
\end{lem}
\proof{} First, we claim that 
\begin{equation} \label{G is small o of F} 
\lim_{u\to+\infty}\frac{G(u)}{F(u)}=0.
\end{equation} 
Indeed, fix $\eps>0$ and recalling that \eqref{KO} holds, choose $M>0$ so large that $\int_{M}^{+\infty}\frac{dt}{\sqrt{2F(t)}}<\eps$. By the definition of $G$, there exists a constant $C_{M}$ such that
$$
G(u) = C_{M} + \int_{M}^{u}{\sqrt{2F(t)}}{dt}.
$$
Since $F$ is nondecreasing it follows that
$$
G(u)\le C_{M} + 2F(u)\int_{M}^{u}\frac{dt}{\sqrt{2F(t)}}\le C_{M} +2\eps F(u).
$$
Dividing by $F(u)$ and letting $u\to+\infty$, \eqref{G is small o of F} follows. 
Next, we claim that
\begin{equation} \label{g is small o of F} 
\lim_{u\to+\infty}\frac{g(u)}{F(u)}=0.
\end{equation}
Note that by \eqref{equation g in u}, $g(u)$ is increasing, thus bounded below by a constant $c$ as $u\to+\infty$. Hence, by \eqref{equation g in u full}, 
$$
\frac{dg}{du}\le \frac{N-1}{r}\sqrt{2(F(u)-c)} \le 2(N-1)\sqrt{2(F(u)-c)},
$$
where the last inequality holds if $r>1/2$ i.e. if $u$ is sufficiently large.
Integrating on a given interval $(u_{0},u)$, we obtain
\begin{equation*}  
c\le g(u) \le g(u_{0}) + 2(N-1)\int_{u_{0}}^{u}\sqrt{2(F(t)-c)}\;dt.
\end{equation*} 
Using \eqref{G is small o of F} and the fact that $\lim_{t\to+\infty}\frac{\sqrt{2F(t)}}{\sqrt{2(F(t)-c)}}=1$, we deduce \eqref{g is small o of F}.  
Now that \eqref{g is small o of F} has been established, we return to \eqref{equation g in u full} and infer that given $\eps>0$, we have for sufficiently large $u$,
$$
\frac{dg}{du}\ge \frac{N-1}{r}\sqrt{2(1-\eps)F(u)} \ge (N-1)\sqrt{2(1-\eps)F(u)}
$$  
and
$$
\frac{dg}{du}\le \frac{N-1}{r}\sqrt{2(1+\eps)F(u)} \le\frac{N-1}{1-\eps}\sqrt{2(1+\eps)F(u)}.
$$
Integrating the above, we finally obtain for large $u$,
$$
(1-\eps)(N-1)\int_{u_{0}}^{u}\sqrt{2F(t)}\;dt \le g(u)-g(u_{0}) \le (1+\eps)^{3/2}(N-1)\int_{u_{0}}^{u}\sqrt{2F(t)}\;dt
$$ 
and \eqref{g is equivalent to G} follows. The fact that $g(u)\to+\infty$ as $u\to+\infty$ follows automatically. 
\hfill\qed
Next, we prove that given two solutions $u_{1},u_{2}$, the corresponding error terms $g_{1}, g_{2}$ given by \eqref{definition g} differ by a bounded quantity.
\begin{lem} \label{lemma two} 
Let $u_{1}, u_{2}$ denote two solutions of \eqref{asymptotic 1}. Introduce $v_{i}=\frac{du_{i}}{dr}$ and 
$$
g_{i}= -\frac{v_{i}^2}{2} + F(u_{i}), \quad\text{for $i=1,2$.}
$$ 
Then, $g_{1}-g_{2}$ is bounded.
\end{lem}
\proof{}
We have seen in the course of the proof of Lemma \ref{lemma one} that each $g_{i}$ (seen as a function of a variable $u$ lying in some interval $(M_{i},+\infty)$)  solves the differential equation \eqref{equation g in u full}. 
%By \eqref{g is equivalent to G} and \eqref{G is small o of F}, $g_{i} = o(F(u))$, so
%\begin{align*}
%\sqrt{2(F(u)-g_{i})} &= \sqrt{2F(u_{i})}\sqrt{1-\frac{g_{i}}{F(u)}}\\
%&= \sqrt{2F(u)}\left(1-\frac12\frac{g_{i}}{F(u)}+ o\left(\frac{g_{i}}{F(u)}\right)\right)\\
%&= \sqrt{2F(u)} -\frac{g_{i}}{\sqrt{2F(u)}} + o\left(\frac{g_{i}}{\sqrt{2F(u)}}\right).
%\end{align*}
%By \eqref{equation g in u full}, we conclude that
%$$
%\frac{dg_{i}}{du} + \frac{N-1}r \frac1{\sqrt{2F(u)}} g_{i}  = \frac{N-1}r \sqrt{2F(u)} + 
%$$ 
Letting $w=g_{1}-g_{2}$, it follows that
\begin{align*}
\frac{dw}{du} &= \frac{N-1}r\left(\sqrt{2(F(u)-g_{1})} - \sqrt{2(F(u)-g_{2})}\right) \\
&= -\frac {2(N-1)}r\frac1{\sqrt{2(F(u)-g_{1})} +\sqrt{2(F(u)-g_{2})}}w.
\end{align*}
So, for some constants $C,u_{0}$,
$$
w = C\exp\left( -\int_{u_{0}}^{u} \frac{2(N-1)}r   \frac1{\sqrt{2(F(t)-g_{1})} +\sqrt{2(F(t)-g_{2})}}\;dt\right),
$$
whence $\left| w \right| \le \left| C \right| $ for $u>u_{0}$.
%By \eqref{g is equivalent to G} and \eqref{G is small o of F}, $g_{i} = o(F(u))$ as $u\to+\infty$. If $u_{0}$ is chosen so large that $r>1/2$ and $g_{i}<\frac12 F(t)$ for $t>u_{0}$, we deduce that
%$$
%\left| w \right| \le C \exp\left( \frac{(N-1)}2\int_{u_{0}}^{u}  \frac1{\sqrt{2F(t)}} \;dt\right)\le C',
%$$
%where we used \eqref{KO} in the last inequality.  
\hfill\qed
\noindent{\bf Completion of the proof of Theorem \ref{theorem one}} 

Let $u_{1}, u_{2}$ denote two solutions of \eqref{asymptotic 1}. By \eqref{equation g in u full}, each $u_{i}$, $i=1,2$, solves
$$ 
\frac{du_{i}/dr}{\sqrt{2(F(u_{i})-g_{i})}} =1.
$$
Integrating, we obtain
$$
\int_{u_{1}}^{+\infty}\frac1{\sqrt{2(F(t)-g_{1})}}\; dt = 1-r = \int_{u_{2}}^{+\infty}\frac1{\sqrt{2(F(t)-g_2)}}\; dt.
$$
Without loss of generality, for a given $r$ we may assume $u_2(r)\ge u_{1}(r)$. Split the left-hand side integral : $\int_{u_{1}}^{+\infty} = \int_{u_{1}}^{u_{2}} + \int_{u_{2}}^{+\infty}$. It follows that
\begin{multline*}
\int_{u_{1}}^{u_{2}}\frac1{\sqrt{2(F(t)-g_{1})}}\; dt = \int_{u_{2}}^{+\infty}\left(\frac1{\sqrt{2(F(t)-g_2)}} - \frac1{\sqrt{2(F(t)-g_1)}}\right)\; dt\\
=\int_{u_{2}}^{+\infty}\frac{\sqrt{2(F(t)-g_1)}-\sqrt{2(F(t)-g_2)}}{\sqrt{2(F(t)-g_1)}{\sqrt{2(F(t)-g_2)}}}\;dt\\
=\int_{u_{2}}^{+\infty}\frac{g_{2}-g_{1}}{\sqrt{2(F(t)-g_1)}\sqrt{2(F(t)-g_2)}\left(\sqrt{2(F(t)-g_1)}+\sqrt{2(F(t)-g_2)}\right) }\;dt
\end{multline*}
Recall that by Lemma \ref{lemma one},  $g_{i}=o(F)$ as $t\to+\infty$. Recall also that $g_{2}-g_{1}$ is bounded. So, for sufficiently large values of $u_{2}$, we deduce
\begin{equation}\label{optimal difference}  
\int_{u_{1}}^{u_{2}}\frac1{\sqrt{2F(t)}}\; dt \le C\int_{u_{2}}^{+\infty}\frac{dt}{F(t)^{3/2}}.
\end{equation} 
Since $F$ is increasing, it follows that
\begin{multline*}
0\le\frac{u_2-u_{1}}{\sqrt{F(u_{2})}}\le C\int_{u_{1}}^{u_{2}}\frac1{\sqrt{2F(t)}}\; dt\le C\int_{u_{2}}^{+\infty}\frac{dt}{F(t)^{3/2}}\\
\le \frac{C}{{\sqrt{F(u_{2})}}}\int_{u_{2}}^{+\infty}\frac{dt}{F(t)}.
\end{multline*}
Hence,
$$
0\le {u_{2} - u_{1}} \le C\int_{u_{2}}^{+\infty}\frac{dt}{F(t)},
$$
as stated in Theorem \ref{theorem one}. It remains to prove \eqref{second estimate}. %Since $F$ is monotone, in virtue of Lemma \ref{lemma three}, we may always assume that $u_{1}$ and $u_{2}$ are radial. 
Without loss of generality, we assume $u_{1}(r)\le u_{2}(r)$ so
\begin{align*}
\int_{u_{1}}^{u_{2}} \frac{dt}{\sqrt{F(t)}} &= \int_{u_{1}}^{+\infty} \frac{dt}{\sqrt{F(t)}} - \int_{u_{2}}^{+\infty} \frac{dt}{\sqrt{F(t)}}\\
&= \int_{u_{2}}^{+\infty} \frac{dt}{\sqrt{F(t-(u_{2}-u_{1}))}} - \int_{u_{2}}^{+\infty} \frac{dt}{\sqrt{F(t)}}\\
&=\int_{u_{2}}^{+\infty} \frac{\sqrt{F(t)} - \sqrt{F(t-(u_2-u_{1}))}}{\sqrt{F(t)F(t-(u_2-u_{1}))}}\;dt\\
&=\int_{u_{2}}^{+\infty} \frac{{F(t)} - {F(t-(u_2-u_{1}))}}{\sqrt{F(t)F(t-(u_2-u_{1}))}\left(\sqrt{F(t)} + \sqrt{F(t-(u_2-u_{1}))}\right)}\;dt\\
&\ge (F(u_{2})-F(u_{1}))\int_{u_{2}}^{+\infty} \frac{dt}{F(t)^{3/2}}.
\end{align*}
Recalling \eqref{optimal difference}, \eqref{second estimate} follows.   
\hfill\qed

\section{Gradient estimates}\label{section 4}
\proof{ of Theorem \ref{theorem finer asymptotics}}
%We want to prove \eqref{third estimate}. 
%Let $\Gamma$ denote the fundamental solution of Laplace's equation and $G_{R}(x,y)$ the Green' s function in the ball of radius $R$ i.e. for $x,y\in B_{R}$, $x\neq y$, 
%$$
%G_{R}(x,y) = \left\{
%\begin{aligned}
%\Gamma(\left| x-y \right| ) - \Gamma\left(\frac{\left| y \right| }{R}\left| x - \frac{R^2}{\left| y \right|^2}y\right| \right),&\qquad y\neq0\\
%\Gamma(\left| x \right| ) - \Gamma(R),&\qquad y=0.
%\end{aligned}
%\right.
%$$
Let $w=u_{1}-u_{2}$ denote the difference of two solutions. Without loss of generality, we may assume that $u_{2}$ is the minimal solution of \eqref{main radial}, so that $u_{1}\ge u_{2}$ and $u_{2}$ is radial.  

{\bf Step 1 : estimate of tangential derivatives}

\noindent We begin by proving that any tangential derivative of $w$ is bounded. Since the problem is invariant under rotation and since $u_{2}$ is radial, we need only show that $\frac{\partial u_{1}}{\partial x_{2}}(r,0,\dots,0)$ remains bounded as $r\to1^-$. Given $x=(x_{1},x_{2},x') \in B$ and $\theta>0$ small,  we denote by $x_{\theta}=(x_{1}\cos\theta - x_{2}\sin\theta, x_{1}\sin\theta + x_{2}\cos\theta,x')$ the image of $x$ under the rotation of angle $\theta$ above the $x_{1}$-axis in the $(x_{1},x_{2})$ plane. By the rotation invariance of the Laplace operator, the function $u_{\theta}$ defined for $x\in B$ by $u_{\theta}(x)=u_{1}(x_{\theta})$, solves \eqref{main radial}. Using \eqref{first estimate} and assumption \eqref{assumption gnn}, we deduce that  $w_{\theta}=u_{1}-u_{\theta}$ solves
\begin{equation}\label{equation w}  
\left\{
\begin{aligned}
\Delta w_{\theta} + K w_{\theta} &= \tilde f(u_{1}) - \tilde f(u_{\theta})&\qquad\text{in $B$,}\\
w_{\theta} &= 0 &\qquad\text{on $\partial B$.}
\end{aligned}
\right.
\end{equation} 
By the Maximum Principle on small domains, there exists $R_{0}\in (0,1)$ such that the operator $L=\Delta + K$ is coercive on $B\setminus B_{R_{0}}$. As a consequence, we claim that there exists a constant $C>0$ such that for all  $x\in B\setminus B_{R_{0}}$,
\begin{equation} \label{estimate w theta} 
\left| w_{\theta}(x) \right| \le C \sup_{\partial B_{R_{0}}} \left| w_{\theta} \right|. 
\end{equation} 
Let indeed $\zeta>0$ denote the solution of
 \begin{equation*}
\left\{
\begin{aligned}
\Delta \zeta + K \zeta &= 0 &\qquad\text{in $B\setminus B_{R_{0}}$}\\
\zeta &= 1 &\qquad\text{on $\partial B_{R_{0}}$}\\
\zeta &= 0 &\qquad\text{on $\partial B$.}
\end{aligned}
\right.
\end{equation*}
We shall prove that $z^\pm:=w_{\theta}-\pm\sup_{\partial B_{R_{0}}}  \left| w_{\theta} \right|\zeta$ are respectively nonpositive and nonnegative, which implies that
\eqref{estimate w theta}  holds for the constant $C= \| \zeta \|_{\infty} $. We work say with $z^+$ and assume by contradiction that the open set $\omega = \{x\in B\setminus B_{R_{0}}\; : \; z^+(x)>0 \}$ is non-empty. Restricting the analysis to a connected component, we have
\begin{equation*}
\left\{
\begin{aligned}
\Delta z^+ + K z^+ &= \tilde f(u_{1}) - \tilde f(u_{\theta})\ge 0\qquad&\text{in $\omega$}\\
z^+ &\le 0\qquad&\text{on $\partial \omega$.}
\end{aligned}
\right.
\end{equation*}
By the Maximum Principle, we conclude that $z^+\le 0$ in $\omega$, a contradiction.   We have thus proved \eqref{estimate w theta}. Since $u_{1}\in C^1(\super{B_{R_{0}}})$, we deduce that for some constant $C>0$ and all $x\in B\setminus B_{R_{0}}$,
$$
\left| w_{\theta}(x) \right| \le C \theta. 
$$  
Applying the above inequality at the point $x=(r,0,\dots,0)$, $r\in(R_{0},1)$ and letting $\theta\to0$, we finally deduce that
$$
\left| \frac{\partial u_{1}}{\partial x_{2}}(r,0,\dots,0)\right| \le C\qquad\text{for all $r\in (R_{0},1)$,}
$$
as desired.

{\bf Step 2 : estimate of the radial derivative}

\noindent It remains to control $\partial w/\partial r$. Fix $R\in(0,1)$. 
Let $G_{R}(x,y)$ denote Green's function in the ball of radius $R$. Then, for $x\in B_{R}$,
\begin{equation} \label{representation} 
\begin{aligned}
w(x) &= \int_{\partial B_{R}} \frac{\partial G_{R}}{\partial \nu_{y}}(x,\cdot)w\;d\sigma &+& \int_{B_{R}}G_{R}(x,\cdot) (f(u_{1})-f(u_{2}))\;dy\\
&=: w_{1}(x) & + & \; w_{2}(x).
\end{aligned}
\end{equation} 
We want to let $R\to1$ in the above identity. To do so, we first observe that $w_{1}$ is harmonic. By the Maximum Principle, $\left| w_{1} \right| \le \| w \|_{L^{\infty}(\partial B_{R})}$. By estimate \eqref{first estimate}, we conclude that $w_{1} \to 0$ as $R\to 1$.  To estimate $w_{2}$, we need the following crucial estimate :
\begin{lem} \label{lemma integral estimate} Assume \eqref{assumption gnn}. Then,
$$\sup_{\theta\in S^{N-1}}\int_{0}^1 \left| f(u_{1}) -  f(u_{2}) \right|(r,\theta) \;dr < +\infty.$$ 
\end{lem}
We shall also need the following elementary estimates.
\begin{lem} \label{potential estimates} There exists a constant $C>0$ such that for all $1/2<r,R< 1$ and all $x,y\in B_{R}$,
\begin{equation}\label{kernel scaling}   
%\left\{
\begin{aligned}
G_{R}(x,y) &= R^{2-N} G_{1}\left(\frac{x}{R},\frac{y}{R}\right)\\
%\frac{\partial G_{R}}{\partial\nu_{y}}(x,y)&=R^{1-N}\frac{\partial G_{1}}{\partial\nu_{y}}\left(\frac{x}{R},\frac{y}{R}\right)
\end{aligned}
%\right.
\end{equation}  
\begin{equation} \label{kernel estimate}
\left\{
\begin{aligned}
\int_{\partial B_{r}}G_{1}(x,\cdot) \;d\sigma&\le 1\\
\int_{\partial B_{r}}\left|  \frac{\partial G_{1}}{\partial \left| x \right| }(x,\cdot) \right| \;d\sigma&\le C\\
\end{aligned}
\right.
\end{equation}  
\end{lem}
We postpone the proofs of the above two lemmas and return to \eqref{representation}. Using polar coordinates,  
\begin{align*}
w_2 (x)  &= \int_{B_{R}}G_{R}(x,\cdot)(f(u_{1})-f(u_{2}))\;dy\\
&=\int_{0}^R \left(\int_{\partial B_{r}}G_{R}(x,\cdot) (f(u_{1})-f(u_{2}))d\sigma\right)\;dr
\end{align*}   
By Lemmas \ref{lemma integral estimate} and \ref{potential estimates}, we may easily pass to the limit in the above expression as $R\to1$, so
$$
w(x) = \int_{B}G_{1}(x,\cdot) (f(u_{1})-f(u_{2}))\;dy\\
$$    
Using again Lemmas \ref{lemma integral estimate} and \ref{potential estimates}, we also have that $w$ is differentiable in the $r= \left| x \right| $ variable and
$$
\frac{\partial w}{\partial r}(x) = \int_{B}\frac{\partial G_{1}}{\partial \left| x \right| }(x,\cdot) (f(u_{1})-f(u_{2}))\;dy
$$   
Using polar coordinates again and Lemmas \ref{lemma integral estimate} and \ref{potential estimates}, we finally obtain
\begin{multline*}
\left| \frac{\partial w}{\partial r} \right|  \le C+\\ +\sup_{r\in(1/2,1)}\left(\int_{\partial B_{r}}\left| \frac{\partial G_{1}}{\partial \left| x \right| }(x,\cdot)  \right| \;d\sigma\right)\sup_{\theta\in S^{N-1}}\left( \int_{0}^1  \left| (f(u_{1})-f(u_{2}) \right|(r,\theta)  \;dr\right)\\
\le C. 
\end{multline*}
It only remains to prove Lemmas \ref{lemma integral estimate} and \ref{potential estimates}.    
\proof{ of Lemma \ref{lemma integral estimate}} 
We first deal with the case where $u_{1}, u_{2}$ are radial and $u_{1}\ge u_{2}$. 
By assumption \eqref{assumption gnn}, we have 
$$
 \int_{0}^1 \left| f(u_{1}) -  f(u_{2}) \right| \;dr\le \int_{0}^1   \left( \tilde f(u_{1}) -  \tilde f(u_{2}) \right)\;dr + K \| u_{1} - u_{2} \|_{L^\infty(B)}.
$$ 
Using \eqref{first estimate}, we see that $u_{1}-u_{2}$ is bounded and so it remains to estimate  $\tilde f(u_{1}) -  \tilde f(u_{2})$.
By \eqref{equation g in u full}, each $u_{i}$, $i=1,2$, solves
$$ 
\frac{du_{i}/dr}{\sqrt{2(F(u_{i})-g_{i})}} =1.
$$
We also know by Lemma \ref{lemma one} that $g_{i}=o(F(u_{i}))$. So, 
$$\lim_{r\to1}\frac{du_{i}/dr}{\sqrt{2F(u_{i})}} = 1.$$
Using this fact, as well as Lemma \ref{lemma two} and \eqref{second estimate}, we obtain for $R\in(1/2,1)$,
\begin{multline*}
\int_{0}^R (\tilde f(u_{1}) - \tilde f(u_{2}))\;dr
\le \int_{0}^R   \left( f(u_{1}) -  f(u_{2}) \right)\;dr + K \| u_{1} - u_{2} \|_{L^\infty(B)}\\
\le C \int_{0}^R(f(u_{1}) - f(u_{2}))\frac{du_{1}/dr}{\sqrt{2F(u_{1})}} \;dr + C\\
\le C \int_{0}^R\left(f(u_{1})\frac{du_{1}/dr}{\sqrt{2F(u_{1})}} - f(u_{2})\frac{du_{2}/dr}{\sqrt{2F(u_{2})}} \right)\;dr+\\
+C\int_{0}^Rf(u_{2})\left( \frac{du_{2}/dr}{\sqrt{2F(u_{2})}} - \frac{du_{1}/dr}{\sqrt{2F(u_{1})}} \right)\;dr + C\\
\le C\left(\sqrt{F(u_{1})}-\sqrt{F(u_{2})}\right)(R)+ C+\\ 
+C\int_{0}^Rf(u_{2})\left( \frac{\sqrt{2(F(u_{2})-g_{2})}}{\sqrt{2F(u_{2})}} - \frac{\sqrt{2(F(u_{1})-g_{1})}}{\sqrt{2F(u_{1})}} \right)\;dr\\
\le C\left(\frac{ F(u_{1})-F(u_{2}) }{ \sqrt{F(u_{1})}+\sqrt{F(u_{2}}) }\right)(R)+ C+\\ 
+C\int_{0}^Rf(u_{2}) \frac{\sqrt{F(u_{2})F(u_{1})-g_{2}F(u_{1})} - \sqrt{F(u_{1})F(u_{2})-g_{1}F(u_{2})} }{\sqrt{F(u_{2})F(u_{1})}} 
\;dr\\
\le \frac{C}{ \sqrt{ F(u_{1}(R)) } }  \| F(u_{1})-F(u_{2})\|_{ L^{\infty} (B)}+ C+\\
+C\int_{0}^R   \frac{f(u_{2})}{ \sqrt{F(u_{2})F(u_{1})}} \frac{ g_{1}F(u_{2}) - g_{2}F(u_{1})}{\sqrt{F(u_{2})F(u_{1})} } 
\;dr\\
\le C + C\int_{0}^R   \frac{f(u_{2})}{ \sqrt{F(u_{2})F(u_{1})}} \frac{ (g_{1}-g_{2})F(u_{2}) + g_{2}(F(u_{2})-F(u_{1}))}{\sqrt{F(u_{2})F(u_{1})} }\; dr \\
\le C + C \| g_{1}-g_{2}\|_{L^\infty(B)}\int_{0}^R \frac{f(u_{2})}{F(u_{2})}\;dr + \\
+C\| F(u_{1})-F(u_{2})\|_{L^\infty(B)}\int_{0}^R \frac{g_{2}}{F(u_{2})}\; dr\\
\le C + C\int_{1/2}^R  \frac{f(u_{2})}{F(u_{2})}\frac{du_{2}/dr}{\sqrt{2F(u_{2})}}\;dr + C\\
\le C + C \left(F^{-1/2}(1/2)-F^{-1/2}(R)\right)\le C.
\end{multline*}
This proves the lemma for radial solutions.
To obtain the estimate in the general case, we may always assume that $u_{2}$ is the minimal solution of \eqref{main radial}, so that $u_{2}\le u_{1}$ and $u_{2}$ is radial. By Lemma \ref{lemma three}, up to replacing $f$ by $\tilde f$ given by Remark \ref{remark KO2}, there exists another radial solution $V$ such that $V\ge u_{1}\ge u_{2}$. Using assumption \eqref{assumption gnn}, we have 
\begin{align*}
 \int_{0}^1 \left| f(u_{1}) -  f(u_{2}) \right| \;dr&\le \int_{0}^1   \left( \tilde f(u_{1}) -  \tilde f(u_{2}) \right)\;dr + K \| u_{1} - u_{2} \|_{L^\infty(B)} \\
 &\le \int_{0}^1 \left( \tilde f(V) -  \tilde f(u_{2}) \right)\;dr + K \| u_{1} - u_{2} \|_{L^\infty(B)} \\
&\le \int_{0}^1 \left(  f(V) -  f(u_{2}) \right)\;dr + 2K \| u_{1} - u_{2} \|_{L^\infty(B)}  
\end{align*}
By \eqref{first estimate}, $u_{1}-u_{2}$ is bounded and the result follows from the radial case.    
\hfill\qed
\proof{ of Lemma \ref{potential estimates}} \eqref{kernel scaling} is standard : write the representation formula \eqref{representation} both in $B_{R}$ and in $B_{1}$, change variables in the $B_{1}$ integral and identify the kernels.  Next, we prove that given any $r\in(0,1)$, $\int_{\partial B_{r}}G_{1}(x,\cdot)d\sigma\le 1$. It suffices to show that for any $\phi\in C_{c}(0,1)$,
\begin{equation} \label{estim} 
\int_{0}^1\phi(r)\left(\int_{\partial B_{r}}G_{1}(x,\cdot)d\sigma\right)\;dr \le  \| \phi \|_{L^1(0,1)}. 
\end{equation} 
By definition of Green's function, the left-hand side of the above inequality is the function $v$ solving
\begin{equation*}  
\left\{
\begin{aligned}
-\Delta v  &= \phi&\qquad\text{in $B$,}\\
v &= 0 &\qquad\text{on $\partial B$.}
\end{aligned}
\right.
\end{equation*} 
The above equation can also be integrated directly :
$$
v'(r)=r^{1-N}\int_{0}^r \phi(t) t^{N-1}\;dt,
$$
 whence $\left| v' \right| \le \| \phi \|_{L^1(0,1)}$ and $\left| v \right| \le \| \phi \|_{L^1(0,1)}$ i.e. \eqref{estim} holds. This proves that $\int_{\partial B_{r}}G_{1}(x,\cdot)d\sigma\le 1$.

% as well as $\left| \int_{\partial B_{r}} \frac{\partial G_{1}}{\partial \left| x \right| }(x,\cdot) \;d\sigma\right| \le 1$. 
 
%Let now $k\ge1$ and let $\varphi_{k}$ denote an eigenfunction associated to $\lambda_{k}$, the $k$-th eigenvalue of the Laplace-Beltrami operator on the sphere $\Delta_{S^{N-1}}$ and $\phi\in C_{c}(0,1)$. Then, the solution of 
%\begin{equation*}  
%\left\{
%\begin{aligned}
%-\Delta v  &= \varphi_{k}\phi&\qquad\text{in $B$,}\\
%v &= 0 &\qquad\text{on $\partial B$.}
%\end{aligned}
%\right.
%\end{equation*}

We turn to the second estimate in \eqref{kernel estimate}. Recall that the Green's function in the unit ball is expressed for $x,y\in B$, $x\neq y$, by
\begin{equation}
G_{1}(x,y) = \Gamma\left(\left(R^2+r^2-2Rr\cos\varphi\right)^{1/2}\right) - \Gamma\left(\left(1+R^2r^2-2Rr\cos\varphi\right)^{1/2}\right),
\end{equation} 
where $R= \left| x \right| $, $r= \left| y \right| $, $\varphi$ is the angle formed by the vectors $x$ and $y$ and $\Gamma$ is the fundamental solution of the Laplace operator. Differentiating with respect to $R$, we obtain for some $C_{N}>0$,
\begin{multline}  \label{green in polar coordinates} 
C_{N}\frac{\partial G_{1}}{\partial \left| x \right|}(x,y) = \\ 
\frac{R-r\cos\varphi}{(R^2+r^2-2Rr\cos\varphi)^{N/2}} - \frac{Rr^2-r\cos\varphi}{(1+R^2r^2-2Rr\cos\varphi)^{N/2}}=\\
\frac{R-r + r(1-\cos\varphi)}{\left((R-r)^2+2Rr(1-\cos\varphi)\right)^{N/2}} - \frac{Rr^2-r+r(1-\cos\varphi)}{\left( (1-Rr)^2 + 2Rr(1-\cos\varphi)\right)^{N/2}}=
A -B\\.
\end{multline} 
We estimate $A$ and leave the reader perform similar calculations for $B$.
Clearly, given $\eps>0$, the expression \eqref{green in polar coordinates}   remains uniformly bounded in the range $1/2<R,r<1$, $\eps<\varphi<2\pi-\eps$. Hence,
$$
\int_{\partial B_{r}}%\frac{\partial G_{1}}{\partial \left| x \right| }(x,\cdot)
\left| A \right| \;d\sigma \le  C_{\eps} + C\int_{\partial B_{r}\cap [0<\varphi<\eps ]} 
%\frac{\partial G_{1}}{\partial \left| x \right| }(x,y) 
\left| A \right| d\sigma.
$$
For $y\in \partial B_{r}\cap [0<\varphi<\eps ]$, let $z=z(y)$ denote the intersection of the line $(Oy)$ and the hyperplane $P$ passing through $x$ and tangent to the hypersphere $\partial B_R$. Then, there exists constants $c_{1},c_{2}>0$ such that for all  $y\in \partial B_{r}\cap [0<\varphi<\eps ]$, 
$$c_{1}(1-\cos\phi)\le \left| z-x \right|^2\le c_{2}(1-\cos\phi).$$ 
Hence, letting $B^{N-1}(x, \rho)\subset P$ denote the $N-1$-dimensional ball of radius $\rho>0$ centered at $x$, we obtain
\begin{multline*}
\int_{\partial B_{r}}\left| A \right| \;d\sigma \le  C\left(1+ \int_{B^{N-1}(x,R\sin\eps)}\frac{\left| R-r \right| +C r \left| z-x \right|^2 }{\left(\left| R-r \right| ^2+c \left| z-x \right|^2  \right)^{N/2} }\;dz\right)\\
\le  C\left(1+ \int_{B^{N-1}(O,R_\eps)}\frac{\left| R-r \right| +C  \left| z\right|^2 }{\left(\left| R-r \right| ^2+c \left| z\right|^2  \right)^{N/2} }\;dz\right)\\
\le  C\left(1+ \int_{B^{N-1}\left(O,\frac{R_\eps}{\left| R-r \right| }\right)}\frac{\left| R-r \right| +C \left| R-r \right|^2 \left| z\right|^2 }{\left| R -r \right|^{N} \left(1+c \left| z\right|^2  \right)^{N/2} }\left| R-r \right|^{N-1} \;dz\right)\\
\le C\left(1 + \int_{\R^{N-1}} \frac1{\left(1+c \left| z\right|^2  \right)^{N/2}}\;dz + \left| R-r \right| \int_{B^{N-1}\left(O,\frac{R_{\eps}}{\left| R-r \right| }\right)}  \left| z \right|^{2-N} \;dz\right)\\
\le C.
\end{multline*}
Working similarly with the $B$ term in \eqref{green in polar coordinates}, we finally obtain the desired estimate \eqref{kernel estimate}. 
\hfill\qed

{%\cb
\section{Asymptotics to all orders}\label{section 5}
This section is devoted to the proof of Theorem \ref{theorem asymptotics}.
Our first task consists in applying the Fixed Point Theorem to the functional $\mathcal N$ defined for $v\in \mathcal B(v_{0},\rho)$, $u\in I$ by 
\begin{equation} \label{functional N} 
[\mathcal N(v)](u) = \sqrt{2\left(F(u)-(N-1)\int_{U_{0}}^{u} \frac{v}{r}\;dt\right)},
\end{equation} 
where $r$ is given by \eqref{r of v}. Let us check first that $\mathcal N\left(\mathcal B(v_{0},\rho)\right)\subset \mathcal B(v_{0},\rho)$. Take $v\in \mathcal B(v_{0},\rho)$. Then,
\begin{equation}\label{estimate r}  
1\ge r\ge 1 -\frac1{1-\rho}\int_{U_0}^{+\infty}\frac1{v_{0}}\;dt=1 -\frac1{1-\rho}\int_{U_0}^{+\infty}\frac1{\sqrt{2F}}\;dt.
\end{equation} 
By \eqref{KO}, it follows that for $\rho<1/4$ and $U_{0}$ sufficiently large, $1\ge r\ge1/2$. Hence,
$$
\left| \int_{U_{0}}^{u}\frac vr\; dt \right| \le C\int_{U_{0}}^{u}\sqrt{2F}\;dt= o(F(u)),
$$
where we used Lemma \ref{lemma one}.  So for $U_{0}$ large and $u\ge U_{0}$, 
$$
\left| \frac{N-1}{F(u)}\int_{U_{0}}^{u}\frac vr\;dt \right| \le \rho.
$$
We deduce that
\begin{equation} \label{N estimate}  
\left| \frac{\mathcal N(v)-v_{0}}{v_{0}}\right| = 1- \sqrt{1 - \frac{N-1}{F(u)}\int_{U_{0}}^{u}\frac vr\;dt}  \le \frac12 \left| \frac{N-1}{F(u)}\int_{U_{0}}^{u}\frac vr dt\right| <\rho.
\end{equation} 
Next, we prove that $\mathcal N$ is contractive. Given $v_{1},v_{2}\in \mathcal B(v_{0},\rho)$, let $r_{1}=r(u,v_{1}), r_{2}=r(u,v_{2})$ (where $r$ is given by \eqref{r of v}). Then, by estimate \eqref{estimate r},  $1/2\le r_1,r_{2}\le 1$ and 
\begin{align*}
\left| \frac{\mathcal N(v_{1})-\mathcal N(v_{2})}{v_{0}}\right| &=
\left| \sqrt{ 1 - \frac{N-1}{F(u)} \int_{U_0}^{u} \frac {v_{1}}{r_{1}}\;dt } -
 \sqrt{1 - \frac{N-1}{F(u)}\int_{U_{0}}^{u}\frac {v_{2}}{r_{2}}\;dt} \right| 
\\
&\le C\frac{N-1}{F(u)}\int_{U_{0}}^{u} \left| \frac{v_{1}}{r_{1}} - \frac{v_{2}}{r_{2}}\right| \;dt
\\
&\le \frac{C}{F(u)}\left(\int_{U_{0}}^{u} \left| v_{1} - v_{2} \right|\;dt + \int_{U_{0}}^{u}v_{0} \left| \frac1{r_{1}} - \frac1{r_{2}} \right|\;dt\right)
\\
&\le \frac{C}{F(u)}\left(\rho\int_{U_{0}}^{u} \sqrt{2F}\;dt + \int_{U_{0}}^{u} \sqrt{2F}\left| {r_{1}} - {r_{2}} \right|\;dt\right)
\\
&\le \frac{C}{F(u)}\left(\rho\int_{U_{0}}^{u} \sqrt{2F}\;dt + \int_{U_{0}}^{u}  \sqrt{2F}\left| \int_{t}^{+\infty} \left(\frac1{v_{1}} - \frac1{v_{2}}\right)\;ds\right|\;dt\right)
\\
&\le \frac{C\rho}{F(u)}\left(\int_{U_{0}}^{u} \sqrt{2F}\;dt + \int_{U_{0}}^{u} \sqrt{2F}\left| \int_{t}^{+\infty} \frac1{\sqrt{2F}}\;ds\right|\;dt\right)\\
&\le \frac{C\rho}{F(u)}\int_{U_{0}}^{u} \sqrt{2F}\;dt.
\end{align*}
Using Lemma \ref{lemma one}, we conclude that $\mathcal N$ is contractive in $\mathcal B(v_{0},\rho)$ if $U_{0}$ was chosen large enough in the first place. We may thus apply the fixed point theorem.
  
So, it only remains to prove \eqref{asymptotics of u}. We first observe that the sequence $(v_{k})$ defined by \eqref{expression vk} is asymptotic i.e. $v_{k+1}(u) = v_{k}(u)(1+o(1))$, as $u\to+\infty$. Since $v_{k+1}=\mathcal N(v_{k})$, it suffices to prove that $\mathcal N(v_{0})-v_{0} = o(v_{0})$ and iterate. By \eqref{N estimate},
$$
\left| \frac{\mathcal N(v_{0})-v_{0}}{v_{0}}\right| \le \frac{C}{F(u)} \int_{U_{0}}^{u}  \sqrt{2F}\;dt
$$ 
and the claim follows by Lemma \ref{lemma one}.  So, the sequence $(v_{k})$ is asymptotic and so must be the sequence $(u_{k})$ defined by \eqref{expression uk}. 
We are now in a position to prove \eqref{asymptotics of u}. By Theorem \ref{theorem two}, we may restrict to the case where $u$ is radially symmetric. Let $v=du/dr$. By \eqref{asymptotic 1}, $v$ solves 
$$
\frac{dv}{dr}+ \frac{N-1}{r}v = f(u)
$$
Use the change of variable $u=u(r)$ to get
$$
v\frac{dv}{du} + \frac{N-1}{r}v = f(u).
$$
Integrating, it follows that for some constant $C$
\begin{equation}\label{with C}  
\frac{v^2}2 = F(u)+C - \int_{U_{0}}^{u}\frac{N-1}{r}v\;dt.
\end{equation} 
Up to replacing $F(u)$ by $\tilde F(u)=F(u)+C$ (which is harmless from the point of view of asymptotics), we may assume $C=0$. So it suffices to prove that $v\in \mathcal B(v_{0},\rho)$ to conclude that $v$ coincides with the unique fixed point of $\mathcal N$, whence \eqref{asymptotics of u} will follow. By \eqref{with C} (with $C=0$), $v\le v_{0}$ and so
\begin{align*}
0\le v_{0}-v &\le \sqrt{2F(u)} -\sqrt{2\left(F(u) -\int_{U_{0}}^u \frac{N-1}{r}v_{0}\;dt\right)}
\\
&\le C\frac{\int_{U_{0}}^u \sqrt{2F}\;dt}{\sqrt{2F(u)}}.
\end{align*}
By Lemma \ref{lemma one}, it follows that
$$
0\le \frac{v_{0}-v}{v_{0}}\le v_{0}-v <\rho
$$ 
and $v\in B(v_{0}, \rho)$ as desired. \hfill\qed

\section{Universal blow-up rate}\label{section 6}
In this section, we prove Theorem \ref{theorem universal}, that is we characterize nonlinearities for which the blow-up rate is universal.
\proof{ of Theorem \ref{theorem universal}  }

{\bf Step 1.} We begin by establishing the theorem when $\Omega=B$ is the unit ball. In light of Theorem \ref{theorem two}, it suffices to prove  \eqref{one term only} for one given solution $u$ of \eqref{main radial}, which we may therefore assume to be radial.     
By \eqref{equation g in u full}, we have after integration that
\begin{equation}\label{first}  
\int_{u}^{+\infty}\frac1{\sqrt{2(F(t)-g)}}\; dt = 1-r.
\end{equation} 
By definition of $u_{0}$, we also have
\begin{equation} \label{second} 
\int_{u_{0}}^{+\infty}\frac1{\sqrt{2F(t)}}\; dt = 1-r.
\end{equation} 
Observe that $u\ge u_{0}$,
split the integral in \eqref{second}  as $\int_{u_0}^{+\infty} = \int_{u_{0}}^{u} + \int_{u}^{+\infty}$ and equate \eqref{first} and \eqref{second}. It follows that
\begin{multline*}
\int_{u_0}^{u}\frac1{\sqrt{2F(t)}}\; dt = \int_{u}^{+\infty}\left(\frac1{\sqrt{2(F(t)-g)}} - \frac1{\sqrt{2F(t)}}\right)\; dt\\
=\int_{u}^{+\infty}\frac{\sqrt{2F(t)}-\sqrt{2(F(t)-g)}}{\sqrt{2F(t)}{\sqrt{2(F(t)-g)}}}\;dt\\
=\int_{u}^{+\infty}\frac{g}{\sqrt{2F(t))}\sqrt{2(F(t)-g)}\left(\sqrt{2F(t)}+\sqrt{2(F(t)-g)}\right) }\;dt
\end{multline*}
Recall that by Lemma \ref{lemma one},  $g=o(F)$ as $t\to+\infty$ and $g(u)\sim (N-1)G(u)=(N-1)\int_{0}^{u}\sqrt{2F}\;dt$. So, for sufficiently large values of $u$, we deduce
\begin{equation}\label{optimal difference bis}  
\int_{u_0}^{u}\frac1{\sqrt{2F(t)}}\; dt \le C\int_{u}^{+\infty}\frac{\int_{0}^{t}\sqrt{2F}\;ds}{(2F(t))^{3/2}}\;dt.
\end{equation} 
Since $F$ is nondecreasing, it follows that
$$
0\le\frac{u-u_0}{\sqrt{2F(u)}}\le \int_{u_0}^{u}\frac1{\sqrt{2F(t)}}\; dt\le C\int_{u}^{+\infty}\frac{\int_{0}^{t}\sqrt{2F}\;ds}{(2F(t))^{3/2}}\;dt.
$$
Hence,
$$
0\le {u - u_0} \le C\sqrt{2F(u)}\int_{u}^{+\infty}\frac{\int_{0}^{t}\sqrt{2F}\;ds}{(2F(t))^{3/2}}\;dt,
$$
and \eqref{one term only} follows from \eqref{assumption one term}.

{\bf Step 2.} Next, we prove that \eqref{one term only} holds for general domains $\Omega$. To this end, we combine a standard approximation argument by inner and outer spheres (see e.g. \cite{lm}) and the comparison technique of \cite{ddgr}. Let $u$ denote a solution of \eqref{main} and take a point $x_{0}\in\partial\Omega$. Let $B\subset\Omega$ denote a ball which is tangent to $\partial\Omega$ at $x_{0}$. Shrink $B$ somewhat by letting $B_{\eps} = (1-\eps)B$, $\eps>0$. Observe that $u\in C(\super{B_{\eps}})$ is a subsolution of 
\begin{equation}\label{main eps}  
 \left\{
 \begin{aligned} 
\Delta U &=f(U) &\quad\text{in $B_{\eps}$,}\\
U &= +\infty &\quad\text{on $\partial B_{\eps}$,}
\end{aligned}
\right. 
 \end{equation} 
By Lemma \ref{lemma KO}, there exists a solution $V_{\eps}$ of \eqref{main eps}, such that $V_{\eps}\ge u$ in $B_{\eps}$. Furthermore, $V_{\eps}$ can be chosen to be the minimal solution of \eqref{main eps} such that $V_{\eps}\ge u$ in $B_{\eps}$. In particular, $V_{\eps}$ is radial and $\eps\to V_{\eps}$ is nondecreasing. In addition, $\eps\to V_{\eps}$ is uniformly bounded on compact subsets of $B$ (working as in the proof of Lemma \ref{lemma KO}), so $V_{\eps}$ converges as $\eps\to0$, to a solution $V$ of \eqref{main radial} such that $V\ge u$ in $B$.  
By Step 1, 
$$\lim_{\substack{x\to x_{0}\\x\in B}} V(x)-u_{0}(1-d_{B}(x))=0,$$ where $d_{B}$ denotes the distance to $\partial B$. Since $V\ge u$ and since the above discussion is valid for any point $x_{0}\in\partial\Omega$, we finally obtain
\begin{equation} \label{limsup} 
\limsup_{x\to\partial\Omega} \left[u(x)-u_{0}(1-d(x))\right] \le0,
\end{equation} 
where $d(x)$ is the distance to $\partial \Omega$. Choose now an exterior ball $B\subset \R^N \setminus\super\Omega$ which is tangent to $\partial\Omega$ at $x_{0}$. For $\eps>0$ small and $R>0$ large, the annulus $A_{\eps}=R B \setminus (1-\eps)B$ contains $\Omega$. Let $U_{\eps}$ denote a large solution on $A_{\eps}$, which we may assume to be minimal, radial and bounded above on $\Omega$ by $u$. Again $U_{\eps}\to U$ as $\eps\to0$ where $U$ is a radial large solution in $A=R B \setminus B\supset\Omega$. Repeating the analysis of Step 1. (which was purely local) for the case of a radial solution defined on an annulus rather than a ball, we easily deduce that
$$\lim_{\substack{x\to x_{0}\\x\in B}} U(x)-u_{0}(1-d_{B}(x))=0.$$
Since $u\ge U$ and since the above discussion is valid for any point $x_{0}\in\partial\Omega$, we obtain
\begin{equation} \label{liminf} 
\liminf_{x\to\partial\Omega} \left[u(x)-u_{0}(1-d(x))\right] \ge0.
\end{equation} 
So, by \eqref{liminf} and \eqref{limsup}, we have that \eqref{one term only} holds in any smoothly bounded domain $\Omega$. 
 
{\bf Step 3.} It only remains to prove that \eqref{one term only} fails when \eqref{assumption one term fails} holds. We use Theorem \ref{theorem asymptotics} to compute the second term in the asymptotic expansion of a solution.       
By \eqref{expression vk}, 
\begin{align*}
v_{1}(u) &= \sqrt{2\left(F(u)-(N-1)\int_{{0}}^{u}\sqrt{2F}\;dt(1+o(1))\right)}\\
&=\sqrt{2F(u)}\left(1 - (N-1)\frac{\int_{{0}}^{u} \sqrt{2F}\;dt}{{2F(u)}}(1+o(1))\right),
\end{align*}
whence
\begin{align*}
\frac1{v_{1}} &= \frac1{\sqrt{2F(u)}}\left(1+(N-1)\frac{\int_{{0}}^{u} \sqrt{2F}\;dt}{{2F(u)}}(1+o(1))\right)\\ 
&=  \frac1{\sqrt{2F(u)}} + (N-1)\frac{\int_{{0}}^{u} \sqrt{2F}\;dt}{{(2F(u))^{3/2}}}(1+o(1)).
\end{align*}

Integrating \eqref{expression uk} for $k=1$, it follows that for $r$ close enough to $1$,
\begin{equation} \label{u1 integrated}
\int_{u_{1}}^{+\infty}\frac{dt}{\sqrt{2F}} + (N-1)(1+o(1))\int_{u_{1}}^{+\infty}\frac{\int_{0}^t\sqrt{2F}\;ds}{(2F)^{3/2}}\;dt = 1-r.
\end{equation}  
Recall \eqref{second}, split the integral in \eqref{second}  as $\int_{u_0}^{+\infty} = \int_{u_{0}}^{u} + \int_{u}^{+\infty}$ and equate \eqref{u1 integrated} and \eqref{second} to get
$$
\int_{u_{0}}^{u_{1}}\frac{dt}{\sqrt{2F}} =  (N-1)(1+o(1))\int_{u_{1}}^{+\infty}\frac{\int_{0}^t\sqrt{2F}\;ds}{(2F)^{3/2}}\;dt.
$$
Since $F$ is nondecreasing, we deduce that
\begin{equation}\label{upper diff}  
\frac{u_{1}-u_{0}}{\sqrt{2F(u_{0})}}\ge (N-1)(1+o(1))\int_{u_{1}}^{+\infty}\frac{\int_{0}^t\sqrt{2F}\;ds}{(2F)^{3/2}}\;dt.
\end{equation} 
Note also that 
\begin{equation} \label{lower diff} 
\int_{u_{0}}^{u_{1}}\frac{\int_{0}^t\sqrt{2F}\;ds}{(2F)^{3/2}}\;dt
\le \int_{u_{0}}^{u_1} \frac{t}{2F}\;dt\le \frac{(u_{1}-u_{0})^2}{4F(u_{0})}.
\end{equation} 
Assume by contradiction that $\lim_{r\to1^-}(u_{1}-u_{0})(r)=0$. Then, \eqref{lower diff} implies that 
$$
\int_{u_{0}}^{u_{1}}\frac{\int_{0}^t\sqrt{2F}\;ds}{(2F)^{3/2}}\;dt = o\left(\frac{u_{1}-u_{0}}{\sqrt{2F(u_{0})}}\right).
$$ 
Using this information in \eqref{upper diff}, we obtain that 
$$
\frac{u_{1}-u_{0}}{\sqrt{2F(u_{0})}}\ge (N-1)(1+o(1))\int_{u_{0}}^{+\infty}\frac{\int_{0}^t\sqrt{2F}\;ds}{(2F)^{3/2}}\;dt.
$$  
But \eqref{assumption one term fails} would then lead us to a contradiction with the assumption $\lim_{r\to1^-}(u_{1}-u_{0})(r)=0$. So we must have   
$$
\liminf_{r\to1^-}\left(u_{1}-u_{0}\right) >0.
$$
%Now, by \eqref{asymptotics of u}, as $x\to\partial B$,
%$$
%u = u_{1}(1+o({1}))
%$$ 
and so \eqref{one term only} fails.\hfill\qed }

\section{The first three singular terms}\label{section 7}
In the previous section, we characterized nonlinearities for which only one term in the expansion is singular. In the present section, we calculate implicitely the next two terms in the expansion. We have not tried to characterize those $f$ for which all remaining terms are nonsingular, but this can certainly be achieved. We leave the tenacious reader try her/his hand at this computational problem.

We begin by calculating the leading asymptotics of $v_{1}$, $v_{2}$. By \eqref{expression vk},  we have 
$$
\frac{v_{1}^2}2 = F - (N-1)\int^{u} {\sqrt {2F}}\;(1+o(1))\;dt.
$$
So,
\begin{multline*}
\frac{v_{1}}{\sqrt 2} = \sqrt{F}\left(1-(N-1)\frac{\int^{u} {\sqrt {2F}}\;dt}{F}(1+o(1))\right)^{1/2}=\\
=\sqrt{F}-\frac{N-1}{2}\frac{\int^{u}\sqrt{2F}\;dt}{\sqrt{F}}(1+o(1)).
\end{multline*}
In other words,
$$
v_{1}= \sqrt{2F}-(N-1)\frac{\int^{u}\sqrt{2F}\;dt}{\sqrt{2F}}(1+o(1)).
$$
To calculate $v_{2}$, we introduce some notation. Given a positive measurable function $v$, set
$$
Pv = \int^{u}v\;dt,\quad Qv=\frac{Pv}{v},\quad Rv = \int^{+\infty}_{u}\frac{dt}{v}\quad\text{ and }Tv=(N-1)PQv +P(vRv).
$$
$v_{1}$ is then expressed by
$$
v_{1} = v_{0} -(N-1)(1+o(1))Qv_{0},
$$
while $v_{2}$ is given by
\begin{multline*}
\frac{v_{2}^2} 2= F - (N-1)\int^{u}\frac{v_{1}}{1-\int_{t}^{+\infty}\frac{ds}{v_{0}}(1+o(1))}\;dt=
\\ = F - (N-1)\int^{u}(v_{0}-(N-1)Qv_{0} + o(Q(v_{0}))(1-Rv_{0}+o(Rv_{0}))\;dt=
\\= F -(N-1)Pv_{0} +(N-1)Tv_{0}(1+ o(1)).
\end{multline*}
So,
\begin{multline*}
v_{2} = \left(2F -2(N-1)Pv_{0} +  2(N-1)T{v_{0}}(1+ o(1))\right)^{1/2}=\\
= v_{0}\left(1-(N-1)\frac{Pv_{0}}{F} +(N-1)\frac{Tv_{0}}{F}(1+o(1))\right)^{1/2}=\\
= v_{0}\left(1-\frac{N-1}{2}\frac{Pv_{0}}{F} + \frac{N-1}{2}\frac{Tv_{0}}{F}-\frac38(N-1)^2\left(\frac{Pv_{0}}{F}\right)^2+\right.\\
+ o (Tv_{0}/F + (Pv_{0}/F)^2)\Big).\\
%=v_{0}-(N-1)Qv_{0}+(N-1)\left(\frac{Tv_{0}}{v_{0}}-\frac34(N-1)\frac{(Pv_{0})^2}{v_{0}^3}\right)(1+o(1)).
\end{multline*}
And so,
\begin{multline*}
\frac1{v_{2}}= \frac1{v_{0}}\left(1+\frac{N-1}{2}\frac{Pv_{0}}{F} - \frac{N-1}{2}\frac{Tv_{0}}{F}+\frac58(N-1)^2\left(\frac{Pv_{0}}{F}\right)^2+\right.\\
+o (Tv_{0}/F + (Pv_{0}/F)^2)\Big)=\\
=\frac1{v_{0}}+(N-1)\frac{Pv_{0}}{v_0^3}+(N-1)\left(-\frac{Tv_{0}}{v_{0}^3}+\frac54(N-1)\frac{(Pv_{0})^2}{v_{0}^5}\right)(1+o(1))=\\
= \frac1{\sqrt{2F}}
+(N-1)\frac{\int^{u}\sqrt{2F}\;dt}{(2F)^{3/2}}+
\\+\frac{(N-1)}{(2F)^{3/2}}
\left(
-\int^{u}\left((N-1)\frac{\int^{t}\sqrt{2F}\;ds}{\sqrt{2F}}+\sqrt{2F}\int^{+\infty}_{u}\frac{ds}{\sqrt{2F}}\right)\;dt+\right.\\
\left.+\frac{5(N-1)}4\frac{\left(\int^{u}\sqrt{2F}\;dt\right)^2}{2F}
\right)(1+o(1))
\end{multline*}
Integrating once more, we finally obtain
\begin{multline*}
1-r=
\int^{+\infty}_{u_{2}(r)}\frac{du}{\sqrt{2F}}+ (N-1)\int^{+\infty}_{u_{2}(r)}\frac{\int^{u}\sqrt{2F}\;dt}{(2F)^{3/2}}\;du+(1+o(1))\times\\
\times(N-1)\int^{+\infty}_{u_{2}(r)}
\left(
-\int^{u}\left((N-1)\frac{\int^{t}\sqrt{2F}\;ds}{\sqrt{2F}}+\sqrt{2F}\int^{+\infty}_{u}\frac{ds}{\sqrt{2F}}\right)\;dt+\right.
\\
\left.+\frac{5(N-1)}4\frac{\left(\int^{u}\sqrt{2F}\;dt\right)^2}{2F}
\right)\;\frac{du}{(2F)^{3/2}}. 
\end{multline*}
This proves Proposition \ref{first three terms}.  
\section{An example: $\mathbf{ f(u)=u^p}$, $\mathbf{ p>1}$}\label{section 8}
Finding the $n$-th term in the expansion for abitrary $n\in\N$ is out of reach for general $f$, simply because of the algorithmic complexity of  calculations. However, when additional information on $f$ is available, one can guess the general form of the expansion and then try to establish it. This is precisely what we do in this section, with the nonlinearity $f(u)=u^p$, $p>1$. 

For notational convenience, we shall work with $F(u)=\frac12 u^{2q}$, where $2q-1=p$, which simply amounts to working with a constant multiple of the original solution.

Recall \eqref{expression vk} and \eqref{expression uk}.  We want to prove inductively that there exists numbers $a_k, b_{k}$ depending on $k,p,N$ only such that
\begin{align}
v_{n} &= u^q\sum_{k=0}^n b_{k}u^{-k(q-1)}+o(u^{q-n(q-1)}),\label{vn} \\
%d&= \sum_{k=1}^{n+1}c_{k}u^{-k(q-1)}+o(u^{-(n+1)(q-1)}),\label{un moins r} \\
u_{n} &= d^{-\frac1{q-1}}\sum_{k=0}^{n}a_{k}d^k + E_{n}(d^{-\frac1{q-1}+n+1}),\label{un} 
\end{align}
where $E_{n}(d^{-\frac1{q-1}+n+1})\sim e_{n} d^{-\frac1{q-1}+n+1}$ for some $e_n\in\R$, as $d\to0^+$.
We have $v_{0}=\sqrt{2F}= u^q$. 
Solving for $u_{0}$ in \eqref{expression uk} yields $u_{0}=c \;d^{-\frac1{q-1}}$. 
%Since $u(r)=u_{0}(r+o(d))$ by \eqref{asymptotics of u}, $d=c u^{-(q-1)}(1+o(1))$. 
So, \eqref{un} and \eqref{vn} hold %\eqref{un moins r} and \eqref{un} 
for $n=0$. Suppose now the result is true for a given $n\in\N$.  In the computations below, the letter $c_{k}$ denotes a number depending on $k,p,N$ only, which value  
may change from line to line. By \eqref{vn}, we have
\begin{multline*}
\frac1{v_{n}}=u^{-q}\left(1+\sum_{k=1}^n b_{k}u^{-k(q-1)}+o(u^{-n(q-1)})\right)^{-1}=
%\nonumber
\\=u^{-q}\left(\sum_{k=0}^n c_{k}u^{-k(q-1)}+o(u^{-n(q-1)})\right). 
\end{multline*}
So,
\begin{align}\label{un sur vn}
\int_{t}^{+\infty}\frac{ds}{v_{n}} &= t^{1-q}\left(\sum_{k=0}^n c_{k}t^{-k(q-1)}\right)+o(t^{-(n+1)(q-1)}) =\\&= \sum_{k=1}^{n+1} c_{k}t^{-k(q-1)}+o(t^{-(n+1)(q-1)}).
\end{align}
It follows that
$$
\frac1{1-\int_{t}^{+\infty}\frac{ds}{v_{n}}}  = \sum_{k=0}^{n+1} c_{k}t^{-k(q-1)}+o(t^{-(n+1)(q-1)}).
$$
Whence,
$$
\frac{v_{n}(t)}{1-\int_{t}^{+\infty}\frac{ds}{v_{n}}} = t^q \sum_{k=0}^{n} c_{k}t^{-k(q-1)}+o(t^{q-n(q-1)}).
$$
And so,
\begin{multline*}
v_{n+1} = \sqrt{2F - (N-1)\int^{u}\frac{v_{n}}{1-\int_{t}^{+\infty}\frac{ds}{v_{n}}}\;dt}=\\
=\sqrt{u^{2q}+u^{q+1} \sum_{k=0}^{n} c_{k}u^{-k(q-1)}+o(u^{1+q-n(q-1)})}=\\
=u^q\left(1+ \sum_{k=1}^{n+1} c_{k}u^{-k(q-1)}+o(u^{-(n+1)(q-1)})\right)^{1/2}=\\
=u^q\sum_{k=0}^{n+1} c_{k}u^{-k(q-1)}+o(u^{q-(n+1)(q-1)}).
\end{multline*}
This proves \eqref{vn}. Integrating \eqref{expression uk}, we obtain
\begin{equation} \label{integratin} 
\int^{+\infty}_{u_{n}}\frac{du}{v_{n}} = d.
\end{equation} 
Now, $v_{n+1} = v_{n} + c_{n+1}u^{q-(n+1)(q-1)}(1+o(1))$. So,
$$
\frac1{v_{n+1}}=\frac1{v_{n}}+c_{n+1} u^{-q-(n+1)(q-1)}(1+o(1)).
$$
It follows that
$$
d=\int_{u_{n+1}}^{+\infty}\frac{du}{v_{n+1}}= \int_{u_{n+1}}^{+\infty}\frac{du}{v_{n}} + 
c_{n+1}u_{n+1}^{-(n+2)(q-1)}(1+o(1)).
$$
In addition, $v_{n}\sim v_{0}$, so $u_{n}\sim u_{0}$, and so $u_{n+1}^{-(q-1)}\sim d$. Using this in the above equation, we get
$$
d+c_{n+1}d^{n+2}(1+o(1)) = \int_{u_{n+1}}^{+\infty}\frac{du}{v_{n}}. 
$$
Recalling that $v_{n}$ is defined by \eqref{integratin} and  satisfies \eqref{un} by induction hypothesis, we conclude that
$$
v_{n+1}=\left(d+c_{n+1}d^{n+2}(1+o(1)\right)^{-\frac1{q-1}}\sum_{k=0}^{n}a_{k}\left(d+c_{n+1}d^{n+2}(1+o(1)\right)^k + E_{n}(d^{-\frac1{q-1}+n+1}).
$$ 
Expanding again the above expression, we finally obtain 
$$
v_{n+1} = d^{-\frac1{q-1}}\sum_{k=0}^{n+1}a_{k}d^k + E_{n+1}(d^{-\frac1{q-1}+n+2}),
$$ 
which proves \eqref{un}. Proposition \ref{utothep} follows. 

\

{\noindent\bf Acknowledgement.} The formulation \eqref{positivity} concerning the positivity of $f$ is due to H. Brezis. The authors gratefully thank him for this suggestion, and for discussing a preliminary version of this work.


\begin{bibdiv}
\begin{biblist}

\bib{bandle}{article}{
   author={Bandle, Catherine},
   title={Asymptotic behavior of large solutions of elliptic equations},
   journal={An. Univ. Craiova Ser. Mat. Inform.},
   volume={32},
   date={2005},
   pages={1--8},
   issn={1223-6934},
   review={\MR{2215890 (2006j:35066)}},
}
	

\bib{bandle-marcus}{article}{
   author={Bandle, Catherine},
   author={Marcus, Moshe},
   title={On second-order effects in the boundary behaviour of large
   solutions of semilinear elliptic problems},
   journal={Differential Integral Equations},
   volume={11},
   date={1998},
   number={1},
   pages={23--34},
   issn={0893-4983},
   review={\MR{1607972 (98m:35049)}},
}

\bib{brezis}{article}{
   author={Brezis, Haim},
   title={Personal communication.}
}

\bib{ddgr}{article}{
   author={Dumont, Serge},
   author={Dupaigne, Louis},
   author={Goubet, Olivier},
   author={R{\u{a}}dulescu, Vicentiu},
   title={Back to the Keller-Osserman condition for boundary blow-up
   solutions},
   journal={Adv. Nonlinear Stud.},
   volume={7},
   date={2007},
   number={2},
   pages={271--298},
   issn={1536-1365},
   review={\MR{2308040 (2008e:35062)}},
}

\bib{gnn}{article}{
   author={Gidas, B.},
   author={Ni, Wei Ming},
   author={Nirenberg, L.},
   title={Symmetry and related properties via the maximum principle},
   journal={Comm. Math. Phys.},
   volume={68},
   date={1979},
   number={3},
   pages={209--243},
   issn={0010-3616},
   review={\MR{544879 (80h:35043)}},
}

\bib{keller}{article}{
   author={Keller, J. B.},
   title={On solutions of $\Delta u=f(u)$},
   journal={Comm. Pure Appl. Math.},
   volume={10},
   date={1957},
   pages={503--510},
   issn={0010-3640},
   review={\MR{0091407 (19,964c)}},
}

\bib{lm}{article}{
   author={Lazer, A. C.},
   author={McKenna, P. J.},
   title={Asymptotic behavior of solutions of boundary blowup problems},
   journal={Differential Integral Equations},
   volume={7},
   date={1994},
   number={3-4},
   pages={1001--1019},
   issn={0893-4983},
   review={\MR{1270115 (95c:35084)}},
}

\bib{osserman}{article}{
   author={Osserman, Robert},
   title={On the inequality $\Delta u\geq f(u)$},
   journal={Pacific J. Math.},
   volume={7},
   date={1957},
   pages={1641--1647},
   issn={0030-8730},
   review={\MR{0098239 (20 \#4701)}},
}

\bib{pv}{article}{
   author={Porretta, Alessio},
   author={V{\'e}ron, Laurent},
   title={Symmetry of large solutions of nonlinear elliptic equations in a
   ball},
   journal={J. Funct. Anal.},
   volume={236},
   date={2006},
   number={2},
   pages={581--591},
   issn={0022-1236},
   review={\MR{2240175 (2007b:35131)}},
}
\end{biblist}
\end{bibdiv}
\end{document}